\documentclass [11pt]{amsart}
\usepackage{multicol}
\usepackage{indentfirst}
\usepackage{latexsym}
\usepackage{bm}
\usepackage{graphicx}
\usepackage{subfigure}
\usepackage[top=1in, bottom=1in, left=1.25in, right=1.25in]{geometry}
\usepackage{epsfig,dsfont,amssymb,amsmath,amsthm,amsfonts,amsbsy,mathrsfs}
\usepackage{psfrag}
\usepackage{subfigure}
\usepackage{epstopdf}
\usepackage{hyperref}
\usepackage{float}
\usepackage{array}
\usepackage{color}
\usepackage{booktabs}

\newtheorem{theorem}{Theorem}[section]

\newtheorem{corollary}[theorem]{Corollary}

\newtheorem{definition}[theorem]{Definition}
\newtheorem{example}[theorem]{Example}

\newtheorem{lemma}[theorem]{Lemma}

\newtheorem{proposition}[theorem]{Proposition}
\newtheorem*{remark}{Remark}

\newtheorem*{outline}{Outline}
\renewcommand{\Im}{\mathrm{Im}}
\newcommand{\Vol}{\mathrm{Vol}}
\renewcommand{\Re}{\mathrm{Re}}

\usepackage[utf8]{inputenc}
\usepackage[english]{babel}
\usepackage{fancyhdr}

\begin{document}
\title[Singular Solutions to the Kapustin-Witten Equations]{Rotationally Invariant Singular Solutions to the Kapustin-Witten Equations}

\author[Siqi He]{Siqi He}
\address{Department of Mathematics, Caltech, 1200 E California Blvd\\ Pasadena, CA 91125}
\email {she@caltech.edu}

\begin{abstract}In the present paper, we find a system of non-linear ODEs that gives rotationally invariant solutions to the Kapustin-Witten equations in 4-dimensional Euclidean space. We explicitly solve these ODEs in some special cases and find decaying rational solutions, which provide solutions to the Kapustin-Witten equations. The imaginary parts of the solutions are singular. By rescaling, we find some limit behavior for these singular solutions. In addition, for any integer $k$, we can construct a 5$|k|$ dimensional family of $C^1$ solutions to the Kapustin-Witten equations on Euclidean space, again with singular imaginary parts. Moreover, we find solutions to the Kapustin-Witten equations over $S^3\times (0,+\infty)$ with the Nahm pole boundary condition.

\end{abstract}

\maketitle

\section{Introduction}
In [\ref{Fivebranes and Knots}], Witten proposed a new physical interpretation of the Jones polynomial and Khovanov homology in terms of counting the solutions of a certain supersymmetric gauge theory in four dimensions and five dimensions.
The BPS equations of N=4 twisted super Yang-Mills theory in four dimensions are called the topological twisted equations [\ref{Electric-Magnetic Duality And The Geometric Langlands Program}] and play an essential role in this framework. The topological twisted equations are:
\begin{equation}
\begin{split}
(F_A-\phi\wedge\phi+\lambda d_A\phi)^+&=0\\
(F_A-\phi\wedge\phi-\lambda^{-1}d_A\phi)^-&=0\\
d^{\star}_A\phi&=0.
\label{TopTwi}
\end{split}
\end{equation}
Witten points out that the most interesting case to study is when $\lambda=-1$. In this case, we obtain the following equations, which we call the Kapustin-Witten equations:
\begin{equation}
\begin{split}
F_A-\phi\wedge\phi-\star d_A\phi&=0\\
d_A^{\star}\phi&=0.
\label{KW}
\end{split}
\end{equation}

In [\ref{Electric-Magnetic Duality And The Geometric Langlands Program}], Kapustin and Witten prove that over a closed manifold, all the regular solutions to the Kapustin-Witten equations are flat $SL(2;\mathbb{C})$ connections. Therefore, the regular solutions to these equations are not so interesting over closed manifolds. However, the Kapustin-Witten equations are interesting over non-compact spaces with singular boundary conditions. Witten's gauge theory approach [\ref{Fivebranes and Knots}] to the Jones polynomial conjectures that the coefficients of the Jones polynomial of a knot are determined by counting the solutions to the Kapustin-Witten equations on $\mathbb{R}^3\times (0,+\infty)$ with the Nahm pole boundary conditions. See also Gaiotto and Witten [\ref{Knot Invariants from Four-Dimensional Gauge Theory}] for an approach to this conjecture. The case of the empty knot is resolved in [\ref{The Nahm Pole Boundary Condition}].

In addition, Taubes studied the compactness properties of the Kapustin-Witten equations [\ref{Compactness theorems for generalizations of the 4-dimensional anti-self dual equations}][\ref{connections on 3-manifolds with bounds on curvature}]. He shows that there can be only two sources of non-compactness. One is the traditional Uhlenbeck bubbling phenomenon [\ref{UhlenbeckCompactness}][\ref{RemoveofSingularities}], and another is the non-compactness coming from the unboundness of the $L^2$ norm of $\phi$.

Therefore, a natural question to ask is whether the Uhlenbeck bubbling phenomenon can appear for solutions to the Kapustin-Witten equations. In addition, do we have a model solution to the Kapustin-Witten equations.

In this paper, we construct some singular solutions to the Kapustin-Witten equations. To be more precise, we prove the following:

\begin{theorem}
After identifying Euclidean space $\mathbb{R}^4$ with the quaternions $\mathbb{H}$, for any real number C, the formulas
\begin{equation}
\left\{
\begin{split}
A(x)=&\Im\Big( \frac{3C}{C^2|x|^4+4C|x|^2+1}\;\bar{x}dx\Big)\\
\phi(x)=&\Im\Big( \frac{3C(C|x|^2+1)}{(C^2|x|^4+4C|x|^2+1)(C|x|^2-1)}\;\bar{x}dx\Big)
\end{split}
\right.
\end{equation}
give solutions to the Kapustin-Witten equations $(\ref{KW})$ with the following properties:

$(1)$ For $C\neq 0$, the solutions are smooth away from $|x|=\frac{1}{\sqrt{C}}$ and decay to 0 when $|x|\rightarrow \infty$.

$(2)$ The solutions have instanton number $0$.

$(3)$ When $C\rightarrow +\infty$, $|F_A|$ converges to a Dirac measure at x=0.

$(4)$ For $C\neq 0$, the pole singularity of $\phi$ at $|x|=\frac{1}{\sqrt{C}}$ cannot be removed by $SU(2)$ gauge transformations.
\end{theorem}

In addition, we also prove the following theorem:
\begin{theorem}
There exists a family of rotationally invariant solutions to the Kapustin-Witten equations on Euclidean $\mathbb{R}^4$ with instanton number $\pm 1$. These solutions are smooth away from a sphere where the real parts are $C^1$ and the imaginary parts are singular.
\end{theorem}

In addition, given an integer $k$, we can generalize the ADHM construction [\ref{ADHMConstructionofinstantons}] and obtain the following theorem:
\begin{theorem}
Given an integer $k$, there exists a $5|k|$ dimensional family of singular solutions to the Kapustin-Witten equations on Euclidean $\mathbb{R}^4$. When $k=\pm 1$, these include the solutions from Theorem $1.2$.
\label{theorem 3}
\end{theorem}
We conjecture that under some non-degeneracy condition, the solutions we obtain in Theorem $\ref{theorem 3}$ have instanton number $k$.

In the last chapter, we observe the relation between the singularity which appears in Thm 1.1 and the Nahm pole boundary condition and get the following theorem.
\begin{theorem}
 There exists two Nahm pole solutions to the Kapustin-Witten equations on $S^3\times (0,+\infty)$, with instanton numbers $\frac{1}{2}$ and $-\frac{1}{2}$.
\end{theorem}

One may ask whether the coefficients of the Jones polynomial of the unknot can be obtained by counting solutions to the Kapustin-Witten equations with the Nahm pole boundary condition on $S^3\times (0,+\infty)$, instead of $\mathbb{R}^3\times(0,+\infty)$. If the solutions found in Theorem 1.4 were unique, this would suggest that a factor of $q^{\frac{1}{2}}+q^{-\frac{1}{2}}$ should be included.

\begin{outline}
In Section $2$, we find a system of non-linear ODEs which will give rotationally invariant solutions to the Kapustin-Witten equations. In Section $3$, we find a first integral of these ODEs and solve them to obtain the solutions in Theorem $1.1$. In Section $4$, we prove the rest part of Theorem $1.1$. In Section $5$, we construct other families of solutions to the Kapustin-Witten equations and prove Theorem $1.2$ and Theorem $1.3$. In Section 6, we build up the relation of our singular solutions and the Nahm pole boundary condition.
\end{outline}

\section{ODEs from the Kapustin-Witten Equations}
\subsection{Background}
In accordance with the philosophy of the ADHM construction [\ref{Geometry of Yang-Mills Field}][\ref{ADHMConstructionofinstantons}] for the anti-self-dual equation, we use quaternions to describe the gauge field in $\mathbb{R}^4$. We begin by briefly recalling the elementary properties of quaternions.

We have three elements $I,\;J,\;K$ satisfying the identities: $I^2=J^2=K^2=-1,\;IJ=-JI=K,\;JK=-KJ=I,\;KI=-IK=J$. A general quaternion $x$ is of the following form:
\begin{equation*}
x=x_1+x_2I+x_3J+x_4K,
\end{equation*}
where $x_1,x_2,x_3,x_4$ are real numbers. After choosing a canonical basis of $\mathbb{R}^4$, we can naturally identify points in $\mathbb{R}^4$ with quaternions. The conjugate quaternion is given by
\begin{equation*}
\bar{x}=x_1-x_2I-x_3J-x_4K
\end{equation*}
and with we have the relation $\overline{xy}=\bar{y}\bar{x}$. In addition, we also know that $x\bar{x}=\bar{x}x=|x|^2=\sum x_i^2$. For $
x=x_1+x_2I+x_3J+x_4K$, the imaginary part of x is $\Im(x):=x_2I+x_3J+x_4K$. Therefore, the Lie group $SU(2)$ can be identified with the unitary quaternions and the Lie algebra $\textbf{su}(2)$ can be identified with the imaginary part of the quaternions.

Using the well known isomorphism of the Lie group $SO(4)$ with $SU(2)\times SU(2)/\sim$, the action of $SO(4)$ on a quaternion $x$ is given by $x\rightarrow axb$, where $a$, $b$ are unitary quaternions.

In order to find rotationally invariant solutions, we assume that the gauge fields of ($\ref{TopTwi}$) ($\ref{KW}$) have the following form:
\begin{equation}
\begin{split}
A(x):&=\Im(f(t)\;\bar{x}dx)\\
\phi(x):&=\Im(g(t)\;\bar{x}dx)\\
t:&=|x|^2.
\label{SolutionType}
\end{split}
\end{equation}

Here $f(t),\;g(t)$ are real functions with variable $t=|x|^2$. Obviously, $t\ge 0$.
\begin{remark}
In the remaining part of the paper, we use $f',g'$ to simplify writing $\frac{df(t)}{dt}$ and $\frac{dg(t)}{dt}$.
\end{remark}
\begin{proposition}
$A(x)$ and $\phi(x)$ defined as in $(\ref{SolutionType})$ are rotationally invariant up to gauge equivalence.
\end{proposition}
\proof It is easy to see that for $a,\;b$ are two unitary quaternions, under the change $x\rightarrow axb$, we obtain $|axb|^2=|x|^2$, $A(axb)=\Im(f(t)\;\overline{axb}\;d(axb))=\bar{b}\Im(f(t)\;\bar{x}dx)b$. Therefore, $A(axb)$ is gauge equivalent to $A(x)$ by a constant gauge transformation. Similarly, we can show $\phi(x)$ is also rotationally invariant up to the same gauge transformation.\qed

\subsection{Basic Properties of Rotationally Invariant Connections}

As the equations $(\ref{KW})$ depend on the metric, we need to be explicit about the metric we choose.

\begin{definition}
A metric $g$ on $\mathbb{R}^4$ is called rotationally invariant if in quaternion coordinate $g=h(t)dx\otimes d\bar{x}$. $h(t)$ here is a positive function, $t=|x|^2$.
\end{definition}

\begin{example}
The Euclidean metric $dx\otimes d\bar{x}$ and the round metric $\frac{4}{(1+t)^2}dx\otimes d\bar{x}$ on $\mathbb{R}^4$ are both rotationally invariant metrics.
\end{example}

\begin{remark}
In the rest of the paper, all the metrics we considered are rotationally invariant.
\end{remark}

Now, we will introduce some basic properties of connections in ($\ref{SolutionType}$).
\begin{lemma}
$\Im(\bar{x}dx\wedge \bar{x}dx)=-\frac{1}{2}|x|^2\;d\bar{x}\wedge dx-\frac{1}{2}\;\bar{x}dx\wedge d\bar{x}x.$
\label{lemma22}
\end{lemma}
\proof Since the wedge product of a real form with itself is zero, we know that $$\Re(\bar{x}dx)\wedge \Re(\bar{x}dx)=0.$$ Since $$\Re(\bar{x}dx)=\frac{\bar{x}dx+d\bar{x}x}{2},$$ we obtain

\begin{equation*}
\begin{split}
0=&\Re(\bar{x}dx)\wedge \Re(\bar{x}dx)\\
=&\frac{(\bar{x}dx+d\bar{x}x)\wedge(\bar{x}dx+d\bar{x}x)}{4}\\
=&\frac{\bar{x}dx\wedge\bar{x}dx+\bar{x}dx\wedge d\bar{x}x+t\;d\bar{x}\wedge dx+d\bar{x}x\wedge d\bar{x}x}{4}.
\end{split}
\end{equation*}


In addition, we have $$\Im(\bar{x}dx\wedge \bar{x}dx)=\frac{\bar{x}dx\wedge \bar{x}dx+d\bar{x}x\wedge d\bar{x}x}{2}.$$ The plus sign on the right hand side of the above identity is because given two quaternion one forms $\omega_1,\omega_2$, we have $\overline{\omega_1\wedge\omega_2}=-\bar{\omega}_2\wedge \bar{\omega}_1$.

The result follows immediately.
\qed

\begin{lemma}
$\Im(\bar{x}dx)\wedge \Im(\bar{x}dx)=\Im(\bar{x}dx\wedge \bar{x}dx)$.
\label{barxdxproduct}
\end{lemma}
\proof
We calculate that
\begin{equation*}
\begin{split}
&\Im(\bar{x}dx)\wedge \Im(\bar{x}dx)\\
=&\frac{(\bar{x}dx-d\bar{x}x)\wedge(\bar{x}dx-d\bar{x}x)}{4}\\
=&\frac{\bar{x}dx\wedge \bar{x}dx-\bar{x}dx\wedge d\bar{x}x-t\;d\bar{x}\wedge dx+d\bar{x}x\wedge d\bar{x}x}{4}\\
=&\frac{\bar{x}dx\wedge \bar{x}dx+d\bar{x}x\wedge d\bar{x}x}{2} \text{    (by Lemma \ref{lemma22})}\\
=&\Im(\bar{x}dx\wedge \bar{x}dx).
\end{split}
\end{equation*}
\qed

\begin{lemma}
$$dx\wedge d\bar{x}=-2((dx_1\wedge dx_2+dx_3\wedge dx_4)I+(dx_1\wedge dx_3+dx_4\wedge dx_2)J+(dx_1\wedge dx_4+dx_2\wedge dx_3)K)$$
$$d\bar{x}\wedge dx=2((dx_1\wedge dx_2+dx_4\wedge dx_3)I+(dx_1\wedge dx_3+dx_2\wedge dx_4)J+(dx_1\wedge dx_4+dx_3\wedge dx_2)K).$$
\label{dx bar dx SD ASD}
\end{lemma}
\proof
By direct computation.
\qed
\begin{remark}
It is easy to see that $dx\wedge d\bar{x}\in\Omega^{2+}$ and $d\bar{x}\wedge dx\in\Omega^{2-}$. Also $\Im(dx\wedge d\bar{x})=dx\wedge d\bar{x}$ and $\Im(d\bar{x}\wedge dx)=d\bar{x}\wedge dx$.
\end{remark}
\subsection{Separating Terms in the Topological-Twisted equations}
Since the equations in (\ref{TopTwi}) are separated into the self-dual parts and the anti-self-dual parts, we also want to separate our calculation into the self-dual parts and the anti-self-dual parts.

\begin{lemma}
For $A(x)$ defined as in $(\ref{SolutionType})$, we have
$$F_A^+=-\frac{1}{2}(f'+f^2)\;\bar{x}dx\wedge d\bar{x}x$$
$$F_A^-=(\frac{1}{2}tf'-\frac{1}{2}tf^2+f)\;d\bar{x}\wedge dx.$$
\label{curvaturePM}
\end{lemma}
\proof
We calculate that
\begin{equation*}
\begin{split}
F_A&=dA+A\wedge A\\
&=d\Im(f\;\bar{x}dx)+\Im(f^2\;\bar{x}dx\wedge \bar{x}dx)\\
&=\Im(df\;\bar{x}dx)+\Im(f\;d\bar{x}\wedge dx)+\Im(f^2\;\bar{x}dx\wedge \bar{x}dx)\\
&=\Im((f'+f^2)\;\bar{x}dx\wedge \bar{x}dx)+\Im((f't+f)\;d\bar{x}\wedge dx)\;(\text{by }x\bar{x}=|x|^2=t)\\
&=-\frac{1}{2}(f'+f^2)t\;d\bar{x}\wedge dx+(f't+f)\;d\bar{x}\wedge dx-\frac{1}{2}(f'+f^2)\;\bar{x}dx\wedge d\bar{x}x \text{    (by Lemma \ref{lemma22})}\\
&=(\frac{1}{2}tf'-\frac{1}{2}tf^2+f)\;d\bar{x}\wedge dx-\frac{1}{2}(f'+f^2)\;\bar{x}dx\wedge d\bar{x}x.
\end{split}
\end{equation*}
The result follows immediately.
\qed

\begin{lemma}
For $\phi(x)$ defined as in $(\ref{SolutionType})$, we have
$$(\phi\wedge\phi)^+=-\frac{1}{2}g^2\;\bar{x}dx\wedge d\bar{x}x$$
$$(\phi\wedge\phi)^-=-\frac{1}{2}g^2t\;d\bar{x}\wedge dx.$$
\label{phi wedge phi PM}
\end{lemma}
\proof
We calculate that
\begin{equation*}
\begin{split}
\phi\wedge \phi=&\Im(g\;\bar{x}dx)\wedge \Im(g\;\bar{x}dx)\\
=&\Im(g^2\;\bar{x}dx\wedge \bar{x}dx)\;(\text{by Lemma } \ref{barxdxproduct})\\
=&-\frac{1}{2}g^2t\;d\bar{x}\wedge dx-\frac{1}{2}g^2\;\bar{x}dx\wedge d\bar{x}x.\;(\text{by Lemma } \ref{lemma22})
\end{split}
\end{equation*}
\qed
\begin{lemma}
For $(A(x),\;\phi(x))$ defined as in $(\ref{SolutionType})$, we have
$$(d_A\phi)^+=-\frac{1}{2}(g'+2fg)\;\bar{x}dx\wedge d\bar{x}x$$
$$(d_A\phi)^-=(\frac{1}{2}g't+g-fgt)\;d\bar{x}\wedge dx.$$
\label{dAplusminus}
\end{lemma}
\proof
We calculate that
\begin{equation*}
\begin{split}
d_A\phi=&d\phi+A\wedge\phi+\phi\wedge A\\
=&d\Im(g\;\bar{x}dx)+\Im(2fg\;\bar{x}dx\wedge \bar{x}dx)\;(\text{by Lemma } \ref{barxdxproduct})\\
=&\Im(dg\;\bar{x}dx)+\Im(g\;d\bar{x}\wedge dx)+\Im(2fg\bar{x}dx\wedge \bar{x}dx)\\
=&\Im((g'+2fg)\;\bar{x}dx\wedge \bar{x}dx)+\Im((tg'+g)\;d\bar{x}\wedge dx)\\
=&(-\frac{1}{2}(g'+2fg)t+(g't+g))\;d\bar{x}\wedge dx-\frac{1}{2}(g'+2fg)\;\bar{x}dx\wedge d\bar{x}x\\
=&(\frac{1}{2}g't+g-fgt)\;d\bar{x}\wedge dx-\frac{1}{2}(g'+2fg)\;\bar{x}dx\wedge d\bar{x}x.
\end{split}
\end{equation*}
\qed

Now, we will discuss the third equation of (\ref{TopTwi}).

At first, we have the following identity:
\begin{lemma}
$x_1\Im(\bar{x})+x_2\Im(\bar{x}I)+x_3\Im(\bar{x}J)+x_4\Im(\bar{x}K)=0.$
\label{Identity2}
\end{lemma}
\proof
We calculate that
\begin{equation*}
\begin{split}
&x_1\Im(\bar{x})+x_2\Im(\bar{x}I)+x_3\Im(\bar{x}J)+x_4\Im(\bar{x}K)\\
=&x_1(-x_2I-x_3J-x_4K)+x_2(x_1I+x_3K-x_4J)+x_3(x_1J-x_2K+x_4I)+x_4(x_1K+x_2J-x_3I)\\
=&0.
\end{split}
\end{equation*}
\qed

\begin{lemma}
For the Hodge star operator correspond to a rotational invariant metric $h(t)dx\otimes d\bar{x}$, we have $d(\Im(\bar{x})\star dx_1+\Im(\bar{x}I)\star dx_2+\Im(\bar{x}J)\star dx_3+\Im(\bar{x}K)\star dx_4)=0$.
\label{Identity1}
\end{lemma}
\proof
By definition, we have $\Im(\bar{x})=-x_2I-x_3J-x_4K$, $\star(dx_1)=hdx_2\wedge dx_3 \wedge dx_4$, therefore $(d\Im(\bar{x}))\star dx_1=0$. Similarly, we have $(d \Im(\bar{x}I))\star dx_2=(d \Im(\bar{x}J))\star dx_3=(d \Im(\bar{x}K))\star dx_4=0$.

Therefore,
\begin{equation}
\begin{split}
&d(\Im(\bar{x})\star dx_1+\Im(\bar{x}I)\star dx_2+\Im(\bar{x}J)\star dx_3+\Im(\bar{x}K)\star dx_4)\\
=&2h^{'}(x_1\Im(\bar{x})+x_2\Im(\bar{x}I)+x_3\Im(\bar{x}J)+x_4\Im(\bar{x}K))\\
=&0.
\end{split}
\end{equation}

\qed

\begin{lemma}
For $(A(x),\;\phi(x))$ defined as in $(\ref{SolutionType})$, we have $A\wedge\star\phi+\star\phi\wedge A=0.$
\label{Aphi+phiA=0}
\end{lemma}
\proof
For $A\wedge \star\phi$, we calculate that
\begin{equation*}
\begin{split}
&A\wedge \star\phi\\
=&fg\;\Im(\bar{x}dx)\wedge\star \Im(\bar{x}dx)\\
=&fg\;(\Im(\bar{x})dx_1+\Im(\bar{x}I)dx_2+\Im(\bar{x}J)dx_3+\Im(\bar{x}K)dx_4)\wedge\\
&(\Im(\bar{x})\star dx_1+\Im(\bar{x}I)\star dx_2+\Im(\bar{x}J)\star dx_3+\Im(\bar{x}K)\star dx_4)\\
=&fg\;(\Im(\bar{x})\Im(\bar{x})dx_1\wedge\star dx_1+\Im(\bar{x}I)\Im(\bar{x}I)dx_2\wedge\star dx_2\\
&+\Im(\bar{x}J)\Im(\bar{x}J)dx_3\wedge\star dx_3+\Im(\bar{x}K)\Im(\bar{x}K)dx_4\wedge\star dx_4).
\end{split}
\end{equation*}
In addition, we calculate that
\begin{equation*}
\begin{split}
&\star\phi\wedge A\\
=&fg\;\star \Im(\bar{x}dx)\wedge \Im(\bar{x}dx)\\
=&fg\;(\Im(\bar{x})\star dx_1+\Im(\bar{x}I)\star dx_2+\Im(\bar{x}J)\star dx_3+\Im(\bar{x}K)\star dx_4)\wedge\\
&(\Im(\bar{x})dx_1+\Im(\bar{x}I)dx_2+\Im(\bar{x}J)dx_3+\Im(\bar{x}K)dx_4)\\
=&fg\;(\Im(\bar{x})\Im(\bar{x})\star dx_1\wedge dx_1+\Im(\bar{x}I)\Im(\bar{x}I)\star dx_2\wedge dx_2\\
&+\Im(\bar{x}J)\Im(\bar{x}J)\star dx_3\wedge dx_3+\Im(\bar{x}K)\Im(\bar{x}K)\star dx_4\wedge dx_4)\\
=&-A\wedge \star\phi.
\end{split}
\end{equation*}
Therefore, we obtain $A\wedge \star\phi+\phi\wedge \star A=0$.
\qed

\begin{proposition}
For $(A(x),\;\phi(x))$ defined as in $(\ref{SolutionType})$, for a Hodge star operator correspond to a rotational invariant metric $h(t)dx\otimes d\bar{x}$, we have $d_A\star \phi=0$.
\label{3rdequationvanish}
\end{proposition}
\proof By definition,
$$d_A\star \phi=d\star\phi+A\wedge\star\phi+\star\phi\wedge A.$$
First, we compute $d\star \phi=0$.

Take $\star_E$ to be the Hodge star operator correspond to the Euclidean metric in $\mathbb{R}^4$, then $\star dx_i=h(t)^2\star_E dx_i$.

By (\ref{SolutionType}), we have
\begin{equation*}
\begin{split}
\phi&=g\;\Im(\bar{x}dx)\\
&=g\;(\Im(\bar{x})dx_1+\Im(\bar{x}I)dx_2+\Im(\bar{x}J)dx_3+\Im(\bar{x}K)dx_4).
\end{split}
\end{equation*}
Therefore, we calculate
\begin{equation*}
\begin{split}
&d\star\phi\\
&=d(gh^2)\;((\Im(\bar{x})\star_E dx_1+\Im(\bar{x}I)\star_E dx_2+\Im(\bar{x}J)\star_E dx_3+\Im(\bar{x}K)\star_E dx_4))\\
&+g\;d(\Im(\bar{x})\star_E dx_1+\Im(\bar{x}I)\star_E dx_2+\Im(\bar{x}J)\star_E dx_3+\Im(\bar{x}K)\star_E dx_4)\\
&=\frac{\partial (gh^2)}{\partial x_1}\;\Im(\bar{x})dx_1\wedge\star_E dx_1+\frac{\partial (gh^2)}{\partial x_2}\;\Im(\bar{x}I)dx_2\wedge\star_E dx_2\\
&+\frac{\partial (gh^2)}{\partial x_3}\;\Im(\bar{x}J)dx_3\wedge\star_E dx_3+\frac{\partial (gh^2)}{\partial x_4}\;\Im(\bar{x}K)dx_4\wedge\star_E dx_4\;(\text{by Lemma \ref{Identity1}})\\
&=2(gh^2)'\;(x_1\Im(\bar{x})+x_2\Im(\bar{x}I)+x_3\Im(\bar{x}J)+x_4\Im(\bar{x}K))dx_1\wedge dx_2\wedge dx_3\wedge dx_4\;(\text{by Lemma \ref{Identity2}})\\
&=0.
\end{split}
\end{equation*}

Combining this with Lemma \ref{Aphi+phiA=0}, we finish the proof.

\qed


\subsection{ODEs from the Kapustin-Witten Equations}
Recall that the topological twisted equations (\ref{TopTwi}) are equivalent to the following:
\begin{equation*}
\begin{split}
F_A^+-(\phi\wedge \phi)^+&=-\lambda (d_A\phi)^+\\
F_A^--(\phi\wedge \phi)^-&=\lambda^{-1}(d_A\phi)^-\\
d^{\star}_A\phi&=0.
\end{split}
\end{equation*}

By Proposition \ref{3rdequationvanish}, we know that $d^{\star}_A\phi=0$ is always satisfied under our assumption (\ref{SolutionType}).

Combining Lemma \ref{curvaturePM} and \ref{dAplusminus}, we obtain the following ODEs:
\begin{equation}
\left\{
\begin{split}
&f'+\lambda g'+f^2-g^2+2\lambda fg=0\\
&tf'-t\lambda^{-1}g'+2f-2\lambda^{-1}g+g^2t-f^2t+2tfg\lambda^{-1}=0.
\label{KapustinWitten=1ODEs}
\end{split}
\right.
\end{equation}

To summarize the previous computation, we have the following theorem:
\begin{theorem}
Given a solution $(f(t),\;g(t))$ to the ODEs $(\ref{KapustinWitten=1ODEs})$, taking $A(x)=\Im(f(x)\;\bar{x}dx)$ and $\phi(x)=\Im(g(x)\;\bar{x}dx)$ gives a solution to the topological twisted equations $(\ref{TopTwi})$.
\end{theorem}

By some linear transformations, we obtain the following ODEs:
\begin{equation}
\left\{
\begin{split}
&(\lambda+\lambda^{-1})tf'+2\lambda f-(\lambda-\lambda^{-1})(tf^2-tg^2)-2g+4fgt=0\\
&(\lambda+\lambda^{-1})tg'+2\lambda^{-1} g-2f+(\lambda-\lambda^{-1})2fgt+2t(f^2-g^2)=0.
\label{TwistODEs}
\end{split}
\right.
\end{equation}

Taking $\lambda=-1$, we obtain
\begin{equation}
\left\{
\begin{split}
&tf'+f+g-2fgt=0\\
&tg'+g+f-t(f^2-g^2)=0.
\label{lambda=1originalODEs}
\end{split}
\right.
\end{equation}

We call the equations (\ref{lambda=1originalODEs}) the $\mathbf{Kapustin}$-$\mathbf{Witten}$ $\mathbf{ODEs}$.

\begin{remark}
The equations $(\ref{TwistODEs})$ are degenerate at t=0, which means that we may not have the uniqueness theorem for a given initial value. Given a solution $(f(t),\;g(t))$ to $(\ref{TwistODEs})$, if we assume $(f(t),\;g(t))$ is well-defined near $t=0$, then we can take $t\rightarrow 0$ in both sides of the equations $(\ref{TwistODEs})$ and we obtain that $\lambda f(0)=g(0)$.
\end{remark}
\section{Explicit Solutions for ODEs}
In this section, we will discuss some properties of the equations (\ref{TwistODEs}) and explicitly solve the equations (\ref{TwistODEs}) for some special cases.
\subsection{Change of Variables}

We are going to simplify the equations (\ref{TwistODEs}) by a change of variables.

Take $\tilde{f}(t):=tf(t),\;\tilde{g}(t):=tg(t)$, then the equations (\ref{TwistODEs}) become
\begin{equation}
\left\{
\begin{split}
&\frac{\lambda+\lambda^{-1}}{2}t\tilde{f}'=\frac{\lambda-\lambda^{-1}}{2}(\tilde{f}^2-\tilde{g}^2-\tilde{f})+\tilde{g}-2\tilde{f}\tilde{g}\\
&\frac{\lambda+\lambda^{-1}}{2}t\tilde{g}'=\frac{\lambda-\lambda^{-1}}{2}(\tilde{g}-2\tilde{f}\tilde{g})+(\tilde{f}-\tilde{f}^2+\tilde{g}^2).
\label{fgtODEs}
\end{split}
\right.
\end{equation}

\begin{remark}
It is easy to see that if $(f,g)$ is a solution for some parameter $\lambda_0$, then $(f,-g)$ is a solution for $-\lambda_0$. This is compatible with changing the orientation of the manifold in the topological twisted equation $(\ref{TopTwi})$.
\end{remark}

Taking $u(t):=\tilde{f}(t)-\frac{1}{2}$, $v(t):=\tilde{g}(t)$, we obtain
\begin{equation}
\left\{
\begin{split}
&\frac{\lambda+\lambda^{-1}}{2}tu'=\frac{\lambda-\lambda^{-1}}{2}(u^2-v^2-\frac{1}{4})-2uv\\
&\frac{\lambda+\lambda^{-1}}{2}tv'=-\frac{\lambda-\lambda^{-1}}{2}2uv-(u^2-v^2-\frac{1}{4}).
\label{EZODEs}
\end{split}
\right.
\end{equation}

In order to obtain an autonomous ODE systems, we take $\tilde{u}(s):=u(e^s)$ and $\tilde{v}(s):=v(e^s)$. We obtain
\begin{equation}
\left\{
\begin{split}
&\frac{\lambda+\lambda^{-1}}{2}\tilde{u}'=\frac{\lambda-\lambda^{-1}}{2}(\tilde{u}^2-\tilde{v}^2-\frac{1}{4})-2\tilde{u}\tilde{v}\\
&\frac{\lambda+\lambda^{-1}}{2}\tilde{v}'=-\frac{\lambda-\lambda^{-1}}{2}2\tilde{u}\tilde{v}-(\tilde{u}^2-\tilde{v}^2-\frac{1}{4}).
\label{KapustinWittenparametrizedODEs}
\end{split}
\right.
\end{equation}
Here $\tilde{u}':=\frac{d\tilde{u}(s)}{ds}$ and $\tilde{v}':=\frac{d\tilde{v}(s)}{ds}$.

\subsection{Some Basic Properties}
Even though the equations (\ref{KapustinWittenparametrizedODEs}) are non-linear, we can find a first integral which can simplify the equations in some special cases.
\begin{proposition}
For $(\tilde{u},\tilde{v})$ a solution to the equations $(\ref{KapustinWittenparametrizedODEs})$ , $I(\tilde{u},\tilde{v})=\frac{1}{3}\tilde{u}^3-\tilde{u}\tilde{v}^2-\frac{1}{4}\tilde{u}-\frac{\lambda-\lambda^{-1}}{2}(\tilde{v}^3-\tilde{u}^2\tilde{v}+\frac{1}{4}\tilde{v})$ is a first integral for the equations $(\ref{KapustinWittenparametrizedODEs})$.\label{FirstintegralforKW}
\end{proposition}
\proof
We calculate that
\begin{equation}
\begin{split}
&(\frac{1}{3}\tilde{u}^3-\tilde{u}\tilde{v}^2-\frac{1}{4}\tilde{u})'\\
=&\tilde{u}^2\tilde{u}'-\tilde{u}'\tilde{v}^2-2\tilde{u}\tilde{v}\tilde{v}'-\frac{1}{4}\tilde{u}'\\
=&\tilde{u}'(\tilde{u}^2-\tilde{v}^2-\frac{1}{4})-2\tilde{u}\tilde{v}\tilde{v}'\\
=&\frac{\lambda-\lambda^{-1}}{\lambda+\lambda^{-1}}((\tilde{u}^2-\tilde{v}^2-\frac{1}{4})^2+(2\tilde{u}\tilde{v})^2).
\end{split}
\end{equation}
We calculate that
\begin{equation}
\begin{split}
&(\tilde{v}^3-\tilde{u}^2\tilde{v}+\frac{1}{4}\tilde{v})'\\
=&\tilde{v}^2\tilde{v}'-\tilde{v}'\tilde{u}^2-2\tilde{v}\tilde{u}\tilde{u}'+\frac{1}{4}\tilde{v}'\\
=&-(\tilde{u}^2-\tilde{v}^2-\frac{1}{4})\tilde{v}'-\tilde{u}'(2\tilde{u}\tilde{v})\\
=&\frac{2}{\lambda+\lambda^{-1}}((\tilde{u}^2-\tilde{v}^2-\frac{1}{4})^2+(2\tilde{u}\tilde{v})^2).
\end{split}
\end{equation}
The proposition follows immediately.
\qed

~\\Since we would like our solution to exist near $t=0$, recalling that $\tilde{u}(s)=e^sf(e^s)-\frac{1}{2}$, $\tilde{v}(s)=e^sg(e^s)$, we obtain the following restrictions: $\lim_{s\rightarrow -\infty}\tilde{u}(s)=-\frac{1}{2}$ and $\lim_{s\rightarrow -\infty}\tilde{v}(s)=0$. Therefore, combining this with Proposition \ref{FirstintegralforKW}, we have the following identity:
\begin{equation}
\frac{1}{3}\tilde{u}^3-\tilde{u}\tilde{v}^2-\frac{1}{4}\tilde{u}-\frac{\lambda-\lambda^{-1}}{2}(\tilde{v}^3-\tilde{u}^2\tilde{v}+\frac{1}{4}\tilde{v})=\frac{1}{12}.
\label{firstintegralingeneral}
\end{equation}

With the first integral we have obtained, we can prove the following:

\begin{proposition}
For $\lambda\neq \pm 1$, if $f(t)$ does not blow-up in finite time, then $g(t)$ will not blow-up in finite time.
\end{proposition}
\proof
By the identity (\ref{firstintegralingeneral}), we have
$$\frac{1}{3}\tilde{u}^3-\frac{1}{4}\tilde{u}=\frac{\lambda-\lambda^{-1}}{2}(\tilde{v}^3-\tilde{u}^2\tilde{v}+\frac{1}{4}\tilde{v})+\tilde{u}\tilde{v}^2+\frac{1}{12}.$$
If $f(t)$ does not blow-up in finite time, $\tilde{u}$ will also not blow-up in finite time. If $\tilde{v}$ blows-up in finite time then the right hand side of the identity will be unbounded but the left hand side will be bounded, which gives a contradiction.
\qed

~\\Even though the topological twisted equations (\ref{TopTwi}) are not conformally invariant, we still have that it is invariant under rescaling by a constant, which leads to the following proposition:
\begin{proposition}
If $(f_0(t),g_0(t))$ is a solution of the equations $(\ref{TwistODEs})$, then for any constant $C$, $(Cf_0(Ct),Cg_0(Ct))$ are solutions to the equations $(\ref{TwistODEs})$.
\end{proposition}
\proof By a direct computation.
\qed

\subsection{t'Hooft Solution when $\lambda=0$}
In this subsection, we will prove that we can obtain the t'Hooft solution of Yang-Mills equation from the equations (\ref{TwistODEs}). By taking $\lambda=0$, (\ref{TopTwi}) becomes
\begin{equation}
\begin{split}
(F_A-\phi\wedge\phi)^+&=0\\
(d_A\phi)^-&=0\\
d_A^*\phi&=0.
\end{split}
\end{equation}

If $\phi=0$, then we are just considering the anti-self dual equation $$F^+_A=0.$$

By taking $\lambda=0$, (\ref{TwistODEs}) becomes
\begin{equation}
\left\{
\begin{split}
&f'+f^2-g^2=0\\
&g't+2g-2tfg=0.
\label{SL(2C)ASDODEs}
\end{split}
\right.
\end{equation}

By Theorem 2.12, we know that every solution to the equations ($\ref{SL(2C)ASDODEs}$) will give a solution for the $SL(2;\mathbb{C})$ anti-self-dual equation.

If $g=0$, the equations (\ref{SL(2C)ASDODEs}) have a solution $(f(t),g(t))=(\frac{1}{1+t},0)$. The corresponding gauge fields are ($A(x)=\Im(\frac{1}{1+|x|^2}\bar{x}dx),\phi(x)=0$), which recovers the t'Hooft solution for anti-self-dual equation in [\ref{thooft1}], [\ref{thooft2}].

\begin{remark}
We can also find a solution using the first integral $I(\tilde{u},\tilde{v})=\tilde{v}\tilde{u}^2-\frac{1}{3}\tilde{v}^3-\frac{1}{4}\tilde{v}$. After some computation, we obtain the solution
$(f(t),g(t))=(\frac{t}{t^2-1},\frac{\sqrt{3}}{t^2-1})$ to $(\ref{SL(2C)ASDODEs})$.
\end{remark}

\subsection{Explicit Solutions to the Kapustin-Witten ODEs}
Taking $\lambda=-1$, the equations ($\ref{TwistODEs}$) become
\begin{equation}
\left\{
\begin{split}
&tf'+f+g-2fgt=0\\
&tg'+f+g-t(f^2-g^2)=0.
\end{split}
\right.
\end{equation}

We can find a solution
\begin{equation}
\left\{
\begin{split}
f(t)&=\frac{1}{2t}\\
g(t)&=\frac{\tan(-\frac{1}{2}\ln(t)+C)}{2t}.
\label{badsolution1}
\end{split}
\right.
\end{equation}
However, the solution will have so many poles that 0 will be an accumulation point of singularities, which is not what we want. We hope to find a solution which is well-defined near 0.

From the equations ($\ref{EZODEs}$), we obtain the ODEs corresponding to the Kapustin-Witten equations:
\begin{equation}
\left\{
\begin{split}
tu'&=2uv\\
tv'&=u^2-v^2-\frac{1}{4}.
\label{NonANtoKW}
\end{split}
\right.
\end{equation}


Recalling that $u(t)=tf(t)-\frac{1}{2}$ and $v(t)=tg(t)$, we hope to obtain a solution well-defined near $t=0$. Therefore, we hope to solve (\ref{NonANtoKW}) for the initial value
$(u(0)=-\frac{1}{2},v(0)=0)$.

By taking $\tilde{u}(s):=u(e^{s}),\tilde{v}(s):=v(e^s)$, we obtain an autonomous system of ODEs:
\begin{equation}
\left\{
\begin{split}
&\tilde{u}'=2\tilde{u}\tilde{v}\\
&\tilde{v}'=\tilde{u}^2-\tilde{v}^2-\frac{1}{4}\\
&\lim_{s\rightarrow -\infty}\tilde{u}(s)=-\frac{1}{2}\\
&\lim_{s\rightarrow -\infty}\tilde{v}(s)=0.
\label{ANtoKW}
\end{split}
\right.
\end{equation}

By Proposition $\ref{FirstintegralforKW}$, we the following identity:
\begin{equation}
\tilde{v}^2\tilde{u}-\frac{1}{3}\tilde{u}^3+\frac{1}{4}\tilde{u}=-\frac{1}{12}.
\label{firstintegraltrival}
\end{equation}
Combining ($\ref{ANtoKW}$) and (\ref{firstintegraltrival}), we obtain
\begin{equation}
\left\{
\begin{split}
&12\tilde{v}^2\tilde{u}=(2\tilde{u}+1)^2(\tilde{u}-1)\\
&\tilde{u}'=2\tilde{u}\tilde{v}\\
&\lim_{s\rightarrow -\infty}\tilde{u}(s)=-\frac{1}{2}\\
&\lim_{s\rightarrow -\infty}\tilde{v}(s)=0.
\label{FinalEquationOfKapustinWitten}
\end{split}
\right.
\end{equation}

Assuming $\tilde{u}<0$, take $W(s):=-\tilde{u}(s)$, so we are trying to solve the following ODEs:
$$
W(s)'=\pm\frac{1}{\sqrt{3}}\sqrt{W(s)}\sqrt{W(s)+1}(2W(s)-1).
$$

We first solve $W(s)'=\frac{1}{\sqrt{3}}\sqrt{W(s)}\sqrt{W(s)+1}(2W(s)-1)$,

Taking $H(s):=\frac{1+4W(s)}{2\sqrt{3}\sqrt{W(s)^2+W(s)}}$, we have the following Lemma:
\begin{lemma}
$\frac{1}{1-H^2(s)}\;dH(s)=-ds$
\end{lemma}
\proof
We calculate that
$$H'(s)=\frac{(2W(s)-1)W(s)'}{4\sqrt{3}(W(s)^2+W(s))\sqrt{W(s)^2+W(s)}}.$$
In addition, we calculate that $$1-H(s)^2=-\frac{(2W(s)-1)^2}{12(W(s)^2+W(s))}.$$

Combining this with $W(s)'=\frac{1}{\sqrt{3}}\sqrt{W(s)}\sqrt{W(s)+1}(2W(s)-1)$, the result follows immediately.
\qed

~\\By the previous lemma, $\frac{1}{2}\ln(\frac{H(s)+1}{H(s)-1})=-s+C$. Therefore, $H(s)=\frac{Ce^{-2s}+1}{Ce^{-2s}-1}$. Combining this with $H(s)=\frac{1+4W(s)}{2\sqrt{3}\sqrt{W(s)^2+W(s)}}$,
we find $W(s)=\frac{2-3H^2+3H\sqrt{H^2-1}}{2(3H^2-4)}$.

Therefore, we have $$W(\ln t)=\frac{1}{2}\frac{C^2t^2-2Ct+1}{C^2t^2+4Ct+1}.$$

By definition, $$f(t)=\frac{\frac{1}{2}-W(\ln t)}{t}.$$ We calculate that $$f(t)=\frac{3C}{C^2t^2+4Ct+1}\;(\text{for any constant C}).$$ Putting this into the equations (\ref{ANtoKW}) and taking $g(t)=\frac{\tilde{v}(ln(t))}{t}$, we obtain $$g(t)=\frac{3C(Ct+1)}{(C^2t^2+4Ct+1)(Ct-1)}.$$

For $W(s)'=-\frac{1}{\sqrt{3}}\sqrt{W(s)}\sqrt{W(s)+1}(2W(s)-1)$, we obtain another solution
$(f(t),g(t))=(\frac{1}{t}\frac{C^2t^2+Ct+1}{C^2t^2+4Ct+1},-\frac{3C(Ct-1)}{(C^2t^2+4Ct+1)(Ct+1)}).$

To summarize, by solving the equations (\ref{TwistODEs}) with $\lambda=-1$, we obtain the following proposition:
\begin{proposition}
\begin{equation}
\left\{
\begin{split}
f_1(t)&=\frac{3C}{(Ct)^2+4(Ct)+1}\\
g(t)&=\frac{3C(Ct+1)}{((Ct)^2+4(Ct)+1)(Ct-1)},
\label{solution1ofKapustinWitten}
\end{split}
\right.
\end{equation}

\begin{equation}
\left\{
\begin{split}
f_2(t)&=\frac{1}{t}\frac{C^2t^2+Ct+1}{C^2t^2+4Ct+1}\\
g(t)&=\frac{3C(Ct+1)}{((Ct)^2+4(Ct)+1)(Ct-1)}
\end{split}
\right.
\end{equation}
are two families of solutions to the Kapustin-Witten ODEs $(\ref{lambda=1originalODEs})$.
\label{f(t)1f(t)2}
\end{proposition}

\section{Instanton Number Zero Solutions}
In this section, we will give a complete proof of Theorem 1.1.

\subsection{Computation of Instanton Numbers}
We will now give a formula to compute the instanton number for the rotationally invariant solutions, which will prove property (2) of Theorem 1.1.
\begin{lemma}
$\bar{x}dx\wedge d\bar{x}x\wedge \bar{x}dx\wedge d\bar{x}x=24t^2\;dx_1\wedge dx_2\wedge dx_3\wedge dx_4$,
$d\bar{x}\wedge dx\wedge d\bar{x}\wedge dx=-24\;dx_1\wedge dx_2\wedge dx_3\wedge dx_4.$
\end{lemma}
\proof
By Lemma \ref{dx bar dx SD ASD}, it is just a direct computation.
\qed

Combining this with Lemma $\ref{curvaturePM}$, we obtain that
\begin{equation}
\begin{split}
|F_A^-|^2&=6(tf'+2f-tf^2)^2\\
|F_A^+|^2&=6(f'+f^2)^2t^2.
\label{curvaturenorm}
\end{split}
\end{equation}

Since we are considering the solutions over the non-compact space $\mathbb{R}^4$, the instanton number is defined as:
\begin{definition}
Given that a connection $(A(x),\;\phi(x))$ is a solution to the Kapustin-Witten equations $(\ref{KW})$, the instanton number $k$ for $(A(x),\;\phi(x))$ is $k:=\frac{1}{4\pi^2}\int_{\mathbb{R}^4} tr(F_A\wedge F_A)\in\mathbb{R}$.
\end{definition}

For a rotationally invariant solution as in (\ref{SolutionType}), we have a simple formula to compute the instanton number.
\begin{proposition}
For a globally defined $C^1$ connection $A(x)=\Im(f(t)\;\bar{x}dx)$ over $\mathbb{R}^4$, by taking $\tilde{f}(t):=tf(t)$, the instanton number $k$ satisfies: $$k=6\int_0^{+\infty}\tilde{f}(\tilde{f}-1)\tilde{f}'dt=(2\tilde{f}^3-3\tilde{f}^2)\mid^{+\infty}_{0}.$$
\label{instantonnumber}
\end{proposition}
\proof
We calculate that
\begin{equation*}
\begin{split}
k=&\frac{1}{4\pi^2}\int_{\mathbb{R}^4} tr(F_A\wedge F_A)\\
&=\frac{1}{4\pi^2}\int_{\mathbb{R}^4} (|F_A^+|^2-|F_A^-|^2)\;dx_1\wedge dx_2\wedge dx_3\wedge dx_4\\
&=\frac{1}{4\pi^2}\int_{\mathbb{R}^4} (6(f'+f^2)^2t^2-6(tf'+2f-tf^2)^2)\;dx_1\wedge dx_2\wedge dx_3\wedge dx_4\\
&=\frac{1}{4\pi^2}\int_{\mathbb{R}^4} 24(tf^2-f)(tf'+f)\;dx_1\wedge dx_2\wedge dx_3\wedge dx_4\\
&=\frac{1}{4\pi^2}\int_{\mathbb{R}^4} 24(\frac{1}{t}\tilde{f}(\tilde{f}-1)\tilde{f}')\;dx_1\wedge dx_2\wedge dx_3\wedge dx_4\\
&=\frac{1}{4\pi^2}12\Vol (S^3)\int_0^{+\infty}(\tilde{f}(\tilde{f}-1)\tilde{f}')\;dt\;(\text{by }dx_1\wedge dx_2\wedge dx_3\wedge dx_4=d\;\Vol _{S^3}\frac{1}{2}tdt)\\
&=6\int_0^{+\infty}(\tilde{f}(\tilde{f}-1)\tilde{f}')\;dt\;(\text{since }\Vol (S^3)=2\pi^2)\\
&=(2\tilde{f}^3-3\tilde{f}^2)\mid^{+\infty}_{0}.
\end{split}
\end{equation*}
\qed

\begin{remark}
The previous formula for the instanton number only works for connections with the specific type $A(x)=\Im(f(t)\;\bar{x}dx)$.
\end{remark}
For $A(x)=\Im(f(t)\;xd\bar{x})$, a conjugate form of (\ref{SolutionType}), we have the following corollary:
\begin{corollary}
For a globally defined $C^1$ connection $A=\Im(f(t)\;\bar{x}dx)$ over $\mathbb{R}^4$, by taking $\tilde{f}(t):=tf(t)$, the instanton number $k$ satisfies:
$$k=-6\int_0^{+\infty}\tilde{f}(\tilde{f}-1)\tilde{f}'dt=(3\tilde{f}^2-2\tilde{f}^3)\mid^{+\infty}_{0}.$$
\label{conjugateinstantonnumber}
\end{corollary}
\proof
We can calculate in a similar way and obtain that
\begin{equation}
\begin{split}
|F_A^+|^2&=6(tf'+2f-tf^2)^2\\
|F_A^-|^2&=6(f'+f^2)^2t^2.
\label{curvaturenorm}
\end{split}
\end{equation}

By the same computation as in Proposition \ref{instantonnumber}, we obtain the result.
\qed
\begin{corollary}
The solution $(f_1(t),g(t))=(\frac{3}{t^2+4t+1},\frac{3(t+1)}{(t^2+4t+1)(t-1)})$ to the Kapustin-Witten ODEs $(\ref{lambda=1originalODEs})$ has instanton number zero.

The solution $(f_2(t),g(t))=(\frac{1}{t}\frac{t^2+t+1}{t^2+4t+1},\frac{3(t+1)}{(t^2+4t+1)(t-1)})$ to the Kapustin-Witten ODEs $(\ref{lambda=1originalODEs})$ has instanton number zero.
\end{corollary}
\proof
Defining $\tilde{f}_1(t):=tf_1(t)$, then $\tilde{f}_1(0)=\tilde{f}_1(+\infty)=0$. Therefore, by Proposition $\ref{instantonnumber}$, we know the instanton number of $A(x)=\Im(f_1(t)\;\bar{x}dx)$ is 0.

Similarly, for $\tilde{f}_2(t):=tf_2(t)$, we have $\tilde{f}_2(0)=\tilde{f}_2(+\infty)=1$. Therefore, by Proposition $\ref{instantonnumber}$, we know the instanton number of $A(x)=\Im(f_2(t)\;\bar{x}dx)$ is 0.
\qed

\begin{proposition}
Given that $(\tilde{f}(t),\tilde{g}(t))$ is a solution to the equations $(\ref{fgtODEs})$, if

$(1)$ $\tilde{f}(t)$ is well defined when $t\in [0,+\infty)$.

$(2)$ $\lim_{t\rightarrow 0}\tilde{f}(t)$ and $\lim_{t\rightarrow 0}\tilde{g}(t)$ exist.

$(3)$ $\lim_{t\rightarrow +\infty}\tilde{f}(t)$ and $\lim_{t\rightarrow +\infty}\tilde{g}(t)$ exist.

$(4)$ $\tilde{f}(t)$ exists for $t\in (0,+\infty)$.

Then $A(x)=\Im(\frac{\tilde{f}(t)}{t}\;\bar{x}dx)$ is a connection with instanton number 0, 1 or -1.
\end{proposition}
\proof
After a change of variable and translation, the ODEs (\ref{fgtODEs}) turn into (\ref{KapustinWittenparametrizedODEs}). Equation (\ref{KapustinWittenparametrizedODEs}) is an autonomous system, therefore the limit point must be a equilibrium point of (\ref{KapustinWittenparametrizedODEs}).

There are two equilibrium points, $(\tilde{u},\tilde{v})\in \{(\frac{1}{2},0),(-\frac{1}{2},0)\}$ or equivalently $(\tilde{f},\tilde{g})\in \{(0,0),(1,0)\}$. Therefore, $(\tilde{f}(0),\tilde{g}(0))\in \{(0,0), (1,0)\}$ and $(\tilde{f}(+\infty),\tilde{g}(+\infty))\in \{(0,0), (1,0)\}$. By Proposition $\ref{instantonnumber}$, we know that $A(x)=\Im(\frac{\tilde{f}(t)}{t}\bar{x}dx)$ can only have instanton number 0, 1 or -1.
\qed

\subsection{Bubbling for Instanton Number 0 Singular Solutions to the Kapustin-Witten Equations}
In this subsection, we will prove property (3) of Theorem 1.1: the existence of some bubbling phenomenon for singular solutions.

By previous computation, we know that:
\begin{equation}
\begin{split}
|F_A^-|^2&=6(tf'+2f-tf^2)^2\\
|F_A^+|^2&=6(f'+f^2)^2t^2.
\label{curvaturenorm}
\end{split}
\end{equation}

Consider the solution
\begin{equation}
\left\{
\begin{split}
f(t)&=\frac{3}{t^2+4t+1}\\
g(t)&=\frac{3(t+1)}{(t^2+4t+1)(t-1)}.
\label{solutionKWnoC}
\end{split}
\right.
\end{equation}

The graphs of $f(t)$ and $g(t)$ are depicted here:

\includegraphics[width=16cm, height=5cm]{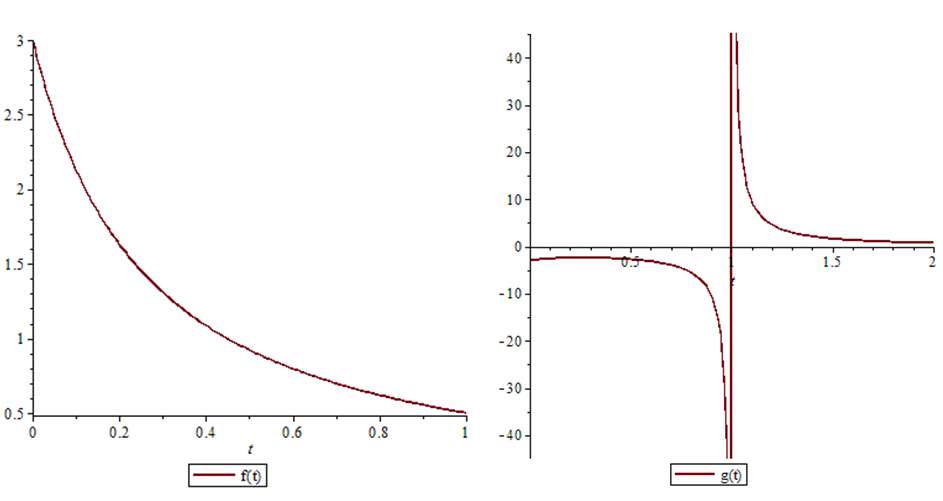}

Combining (\ref{curvaturenorm}) and (\ref{solutionKWnoC}), we obtain that
\begin{equation*}
\begin{split}
|F_A|^2=&|F_A^+|^2+|F_A^-|^2\\
=&(6(tf'+2f-tf^2)^2+6(f'+f^2)^2t^2)\\
=&\frac{108(2t^4+2t^3+t^2+2t+2)}{(t^2+4t+1)^4}.
\end{split}
\end{equation*}

The graph of $|F_A|(t)$ is depicted here:

\includegraphics[width=\textwidth, height=4cm]{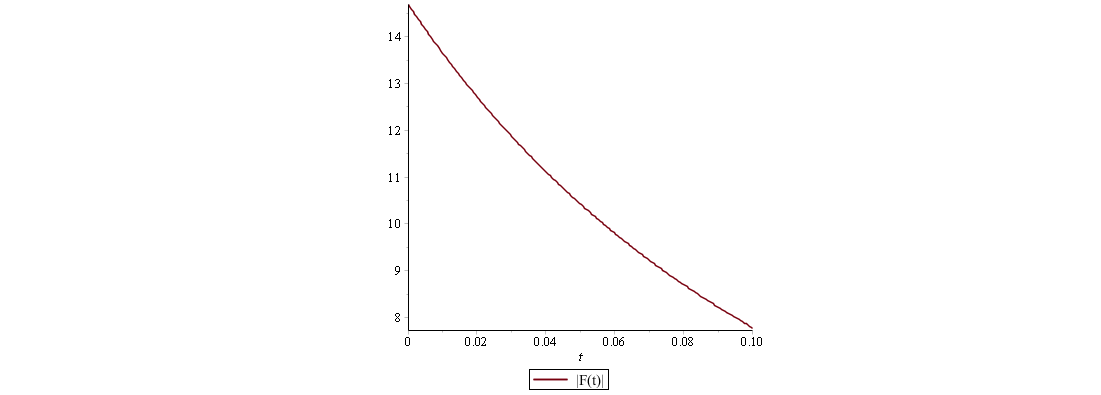}

\begin{proposition}
$|F_A|(t)$ is decreasing and $|F_A|(0)$ is its maximum.
\label{decaymaximal}
\end{proposition}
\proof
A direct computation shows that $\frac{d}{dt}|F_A|^2(t)=-\frac{216(4t^5+5t^4+3t^3+8t^2+19t+15)}{(t^2+4t+1)^5}$, so $|F_A|'<0$ for all $t\geq 0$. Therefore, $|F_A|$ is decreasing and $|F_A|(0)$ is its maximum.
\qed

~\\By Proposition 3.5, for any constant C, we have the solutions
\begin{equation}
\left\{
\begin{split}
f^C(t)&=\frac{3C}{C^2t^2+4Ct+1}\\
g^C(t)&=\frac{3C(Ct+1)}{(C^2t^2+4Ct+1)(Ct-1)}
\label{solutionKW}
\end{split}
\right.
\end{equation}
to the Kapustin-Witten ODEs.

We define $|F_A^C|(t)$ as the curvature norm for $f^C(t)$, then we have the following proposition:

\begin{proposition}
$|F_A^C|(t)=C|F_A|(Ct).$
\end{proposition}
\proof
This is an immediate computation.
\qed

\begin{proposition}
$\lim_{C\rightarrow\infty}|F_A^C|(t)$ is a Dirac measure centered at zero.
\end{proposition}
\proof
By Proposition \ref{decaymaximal}, $|F_A|$ is fast decaying and $|F_A|(0)$ is the maximum of $|F_A|(t)$. Thus by rescaling, we obtain the Dirac measure at $t=0$.
\qed

\subsection{Non-removability of Singularities for $\phi$ by SU(2) Gauge Transformations}
In this subsection, we will prove property (4) of Theorem 1.1. It is suffice to consider the $C=1$ case.

By Lemma \ref{dAplusminus}, we obtain that
$$(d_A\phi)^+=-\frac{1}{2}(g'+2fg)\;\bar{x}dx\wedge d\bar{x}x=\frac{3(t^3+3t+2)}{(t^2+4t+1)(t-1)(t^3+3t^2-3t+1)}\;\bar{x}dx\wedge d\bar{x}x$$
$$(d_A\phi)^-=(\frac{1}{2}g't+g-fgt)\;d\bar{x}\wedge dx=-3\frac{2t^3+3t^2+1}{(t^2+4t+1)(t-1)(t^3+3t^2-3t+1)}d\bar{x}\wedge dx.$$

We calculate that
$$|d_A\phi|^2=\frac{432(2t^8+6t^7+5t^6+2t^5+6t^4+2t^3+5t^2+6t+2)}{(t^2+4t+1)^2(t-1)^2(t^3+3t^2-3t+1)^2}.$$ Therefore, we know $||d_A\phi||_{L^2}$ is unbounded near $t=1$.

Since $||d_A\phi||_{L^2}$ is invariant under the $SU(2)$ gauge action, we know that the singularities of $\phi$ can not be removed by $SU(2)$ gauge transformations.

\section{Non-Zero Instanton Number Solutions}
\subsection{Instanton Number $\pm 1$ Solutions}
In this subsection, we are going to give a proof of Theorem 1.2.

First, we are going to give a construction of an instanton number 1 solution.

By Proposition $\ref{instantonnumber}$, we know that the instanton number is determined by the limit behavior of our connection $A(x)=\Im(f(t)\;\bar{x}dx)$. In order to construct an instanton number $\pm 1$ solution, we only need to construct a solution with different equilibrium points at $t=0$ and $t=+\infty$.

By Proposition \ref{f(t)1f(t)2}, taking $C=1$, we have the following solutions:
\begin{equation}
\left\{
\begin{split}
f_1(t)&=\frac{3}{t^2+4t+1}\\
g(t)&=\frac{3(t+1)}{(t^2+4t+1)(t-1)},
\end{split}
\right.
\end{equation}
\begin{equation}
\left\{
\begin{split}
f_2(t)&=\frac{1}{t}\frac{t^2+t+1}{t^2+4t+1}\\
g(t)&=\frac{3(t+1)}{(t^2+4t+1)(t-1)}.
\end{split}
\right.
\end{equation}

The graphs of $f_1(t)$ and $f_2(t)$ are depicted here:

\includegraphics[width=\textwidth, height=4cm]{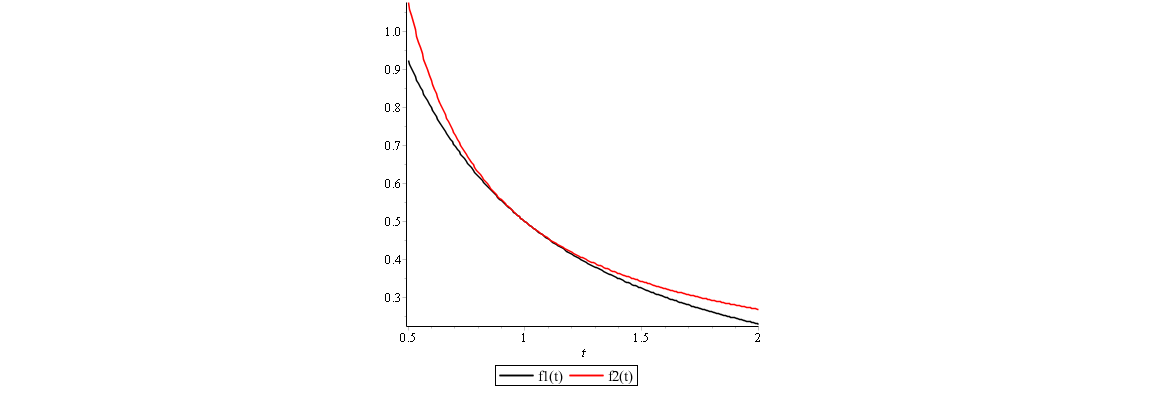}

We hope to glue these two solutions to obtain a new solution.

\begin{proposition}
\begin{equation}
A(x)=\left\{
\begin{split}
&\Im\Big(\frac{3}{t^2+4t+1}\;\bar{x}dx\Big)\;(t\leq 1)\\
&\Im\Big(\frac{1}{t}\frac{t^2+t+1}{t^4+4t+1}\;\bar{x}dx\Big)\;(t\geq 1)
\end{split}
\right.
\end{equation}
\begin{equation}
\begin{split}
\phi(t)=\Im\Big(\frac{3(t+1)}{(t-1)(t^2+4t+1)}\;\bar{x}dx\Big)
\end{split}
\end{equation}
is a solution to the Kapustin-Witten Equations $(\ref{KW})$ which satisfies the following properties:

$(1)$ The solution has instanton number=1.

$(2)$ $A(x)$ is $C^{\infty}$ away from $t=1$ and $C^1$ at $t=1$, $\phi(t)$ is singular at $t=1$.

\label{instantonnumber1solution}
\end{proposition}
\proof

\begin{equation}
f(t)=\left\{
\begin{split}
\frac{3}{t^2+4t+1}\;(t\leq 1)\\
\frac{1}{t}\frac{t^2+t+1}{t^4+4t+1}\;(t\geq 1)
\end{split}
\right.
\end{equation}

Taking $\tilde{f}(t)=tf(t)$, by a direct computation, we know that $\tilde{f}(0)=0$, $\tilde{f}(+\infty)=1$. By Proposition $\ref{instantonnumber}$, we know that the instanton number is equal to 1.

Defining $u(t):=\tilde{f}(t)-\frac{1}{2}$, then
\begin{equation}
u(t)=\left\{
\begin{split}
\frac{1}{2}\frac{t^2-2t+1}{t^2+4t+1}\;(t\leq 1)\\
-\frac{1}{2}\frac{t^2-2t+1}{t^2+4t+1}\;(t\geq 1).
\end{split}
\right.
\end{equation}
By a direct computation, we know that $u(t)$ is a $C^1$ function. Therefore, $A(x)$ is also a $C^1$ connection.
\qed

\begin{remark}
By Corollary $\ref{conjugateinstantonnumber}$, we know
\begin{equation}
A(x)=\left\{
\begin{split}
&\Im\Big(\frac{3}{t^2+4t+1}\;xd\bar{x}\Big)\;(t\leq 1)\\
&\Im\Big(\frac{1}{t}\frac{t^2+t+1}{t^2+4t+1}\;xd\bar{x}\Big)\;(t\geq 1)
\end{split}
\right.
\end{equation}
\begin{equation}
\begin{split}
\phi(x)=\Im\Big(\frac{3(t+1)}{(t-1)(t^2+4t+1)}\;xd\bar{x}\Big)
\end{split}
\end{equation}

is a instanton number $-1$ solution to the Kapustin-Witten equations.
\end{remark}
\subsection{Linear Combination of Solutions}

In this subsection, we aim to generalize the ADHM construction from [\ref{Geometry of Yang-Mills Field}] to obtain higher instanton number solutions to the Kapustin-Witten equations. However, there exists an essential problem to generalizing the instanton number computation method from the anti-self-dual equation case [\ref{Geometry of Yang-Mills Field}]. We conjecture that we will obtain some higher instanton number solutions from this construction.

In view of Corollary \ref{conjugateinstantonnumber}, without loss of generality, we can focus on instanton number $k\geq 0$.

Now, let $\lambda_1,...,\lambda_k$ be $k$ real numbers and $b_1,...,b_k$ be $k$ numbers in $\mathbb{H}$. Take $U:=(\lambda_1(x-b_1),...,\lambda_k(x-b_k))^{T}$ and $U^{\star}$ be the conjugate transpose of $U$. Take $e_0=1,\;e_1=I,\;e_2=J,\;e_3=K$, then for any quaternion $b_i$, we can write $b_i=b_{ij}e_j$.

Now we are going to compute an identity which is parallel to $k=0$ case.

\begin{lemma}
For any $g(t)\in C^1$, $d\star \Im(g(|U|^2)\;U^\star dU)$=0
\label{dstarphi=0 for high instanton}
\end{lemma}
\proof
We calculate that
\begin{equation}
\begin{split}
&d\Im(g(|U|^2)\;\star U^\star dU)\\
=&\sum^k_{i=1}dg(|U|^2)\Im(\lambda_i^2(\bar{x}-\bar{b}_i)\wedge \star d(x-b_i))\;(\text{by Lemma \ref{Identity1}})\\
=&\sum^k_{i=1}\sum^4_{j=1}dg(|U|^2)\lambda_i^2\Im((\bar{x}-\bar{b}_i)e_j)\wedge \star d(x_j-b_{ij})\\
=&\sum^k_{i=1}\lambda_i^2\sum^4_{j=1}\frac{\partial g}{\partial x_j}\;\Im((\bar{x}-\bar{b}_i)e_j)\;dx_j\wedge \star dx_j\\
=&\sum^k_{i=1}\lambda_i^2g ' \sum^4_{j=1} \frac{\partial|U|^2}{\partial x_j}\; \Im((\bar{x}-\bar{b}_i) e_j)\;dx_j\wedge \star dx_j\\
=&\sum^k_{i=1}\lambda_i^2g '\sum^k_{l=1}\lambda_l^2(x_j-b_{lj})\sum_{j=1}^4\Im((\bar{x}-\bar{b}_i)e_j)\;dx_j\wedge \star dx_j\\
=&2g' d\Vol  \sum^k_{i=1}\sum^k_{l=1}\sum^4_{j=1}\lambda_i^2\lambda_l^2(x_j-b_{lj})\Im((\bar{x}-\bar{b}_i)e_j)\\
=&g'd\Vol  \sum^k_{i=1}\sum^k_{l=1}\sum^4_{j=1}\lambda_i^2\lambda_l^2(\sum^4_{j=1}(x_j-b_{lj})\Im((\bar{x}-\bar{b}_i)e_j)+\sum^4_{j=1}(x_j-b_{ij})\Im((\bar{x}-\bar{b}_l)e_j)).
\end{split}
\end{equation}
Therefore, in order to show $d\star \Im(g(|U|^2)\;U^\star dU)$=0, we only need to show $$\sum^4_{j=1}(x_j-b_{lj})\Im((\bar{x}-\bar{b}_i)e_j)+\sum^4_{j=1}(x_j-b_{ij})\Im((\bar{x}-\bar{b}_l)e_j)=0.$$

By translation, without loss of generality, we can assume $b_i=0$. Then we calculate that
\begin{equation}
\begin{split}
&\sum^4_{j=1}(x_j-b_{lj})\Im((\bar{x})e_j)\\
=&\sum^4_{j=1}x_j\Im((\bar{x})e_j)-\sum^4_{j=1}b_{lj}\Im((\bar{x})e_j)\\
=&-\sum^4_{j=1}b_{lj}\Im((\bar{x})e_j).
\end{split}
\end{equation}

For the rest, we calculate that
\begin{equation}
\begin{split}
&\sum^4_{j=1}x_j\Im((\bar{x}-\bar{b}_l)e_j)\\
=&-\sum^4_{j=1}x_j\Im(\bar{b}_le_j)\\
=&-(x_1(-b_{l2}I-b_{l3}J-b_{l4}K)+x_2(b_{l1}I+b_{l3}K-b_{l4}J)\\
&+x_3(b_{l1}J-b_{l2}K+b_{l4}I)+x_4(b_{l1}K+b_{l2}J-b_{l3}I))\\
=&\sum^4_{j=1}b_{lj}\Im((\bar{x})e_j).
\end{split}
\end{equation}

Therefore, we obtain the following identity: $$\sum^4_{j=1}(x_j-b_{lj})\Im((\bar{x}-\bar{b}_i)e_j)+\sum^4_{j=1}(x_j-b_{ij})\Im((\bar{x}-\bar{b}_l)e_j)=0.$$

Combining all the things above, we obtain the lemma.
\qed
\begin{lemma}
For any $f(t),\;g(t)\in C^1$, we have $\Im(f(|U|^2)\;U^\star dU)\wedge \star \Im(g(|U|^2)\;U^\star dU)=0$.
\label{Aphi phiA=0}
\end{lemma}
\proof
Since $f(|U|^2),g(|U|^2)$ are real functions, we only need to show that $\Im(U^\star dU)\wedge \star \Im(U^\star dU)=0$.

We calculate that
\begin{equation}
\begin{split}
&\Im(U^{\star }dU)\wedge \star \Im(U^{\star }dU)\\
=&\sum_{i=1}^k\sum_{j=1}^4\lambda_i^2\Im((\bar{x}-\bar{b}_i)e_j)\;dx_j\wedge \sum_{m=1}^k\sum_{n=1}^4\lambda_m^2\Im((\bar{x}-\bar{b}_m)e_n)\;\star dx_n\\
=&\sum_{i=1}^k\sum_{m=1}^k\sum_{j=1}^4\lambda_i^2\lambda_m^2\Im((\bar{x}-\bar{b}_i)e_j)\Im((\bar{x}-\bar{b}_m)e_j)\;dx_j\wedge \star dx_j\\
=&-\star \Im(U^{\star }dU)\wedge \Im(U^{\star }dU).
\end{split}
\end{equation}
\qed

\begin{corollary}
For any $f(t),g(t)\in C^1$, if $A(x)=\Im(f(t)\;\bar{x}dx)$ and $\phi(x)=\Im(g(t)\;\bar{x}dx)$, then we have $d_A\star \phi=0$.
\label{MultionstatondAstarphi=0}
\end{corollary}
\proof
We have $d_A\star \phi=d\phi+A\wedge \star \phi+\phi\wedge \star A$. This is a direct corollary of Lemma \ref{dstarphi=0 for high instanton} and Lemma \ref{Aphi phiA=0}.
\qed

Taking $f_1(t)=\frac{3}{t^2+4t+1}$, $f_2(t)=\frac{1}{t}\frac{t^2+t+1}{t^2+4t+1}$, $g(t)=\frac{3(t+1)}{(t-1)(t^2+4t+1)}$, we have the following proposition:

\begin{proposition}
\begin{equation}
A(x)=\left\{
\begin{split}
&\Im\Big(f_1(|U|^2)\;U^\star dU\Big)=\Im\Big(\frac{3}{|U|^4+4|U|^2+1}U^\star dU\Big)\;(|U| \leq 1)\\
&\Im\Big(f_2(|U|^2)\;U^\star dU\Big)=\Im\Big(\frac{1}{|U|^2}\frac{|U|^4+|U|^2+1}{|U|^4+4|U|^2+1}U^\star dU\Big)\;(|U| \geq 1)
\end{split}
\right.
\end{equation}
\begin{equation}
\begin{split}
\phi(x)=\Im\Big(g(|U|^2)\;U^\star dU\Big)=\Im\Big(\frac{3(|U|^2+1)}{(|U|^2-1)(|U|^4+4|U|^2+1)}U^\star dU\Big)
\end{split}
\end{equation}
are solutions to the Kapustin-Witten equations.
\label{5DfamilyofSolution}
\end{proposition}
\proof By Corollary \ref{MultionstatondAstarphi=0}, we only need to show that $F_A-\phi\wedge\phi-\star d_A\phi=0$. As before, by replacing $t=|x|^2$ with $|U|^2$, we obtain similar ODEs comparing to ($\ref{KapustinWitten=1ODEs}$).

We obatin that
\begin{equation}
\left\{
\begin{split}
&f(|U|^2)'+\lambda g(|U|^2)'+f(|U|^2)^2-g(|U|^2)^2+2\lambda f(|U|^2)g(|U|^2)=0\\
&tf(|U|^2)'-t\lambda^{-1}g(|U|^2)'+2f(|U|^2)-2\lambda^{-1}g(|U|^2)+g(|U|^2)^2t-f(|U|^2)^2t+2tf(|U|^2)g(|U|^2)\lambda^{-1}=0.
\end{split}
\right.
\end{equation}
The derivative here is the derivative of $|U|^2$. Comparing this with (\ref{TwistODEs}), we are exactly solving the same equations. Therefore, our construction gives solutions to the Kapustin-Witten equations.
\qed
~\\ \begin{proof}[Proof of Theorem 1.3] By our construction, we have the freedom to choose $k$ real numbers $\lambda_1,...,\lambda_k$ and $k$ quaternions $b_1,...,b_k$ in $\mathbb{H}$ in Proposition $\ref{5DfamilyofSolution}$. Therefore, we have a 5$k$ dimension family of solutions to the Kapustin-Witten equations.
\end{proof}

\section{Nahm Pole Boundary Solutions over $S^3\times (0,+\infty)$}
In this section, we will show that our solutions in Section 3 can provide solutions on $S^3\times (0,+\infty)$ with the Nahm pole boundary condition.
\subsection{Nahm Pole Boundary condition}
Now, we will discuss what is the Nahm pole boundary condition to the Kapustin-Witten equations. Denote by $Y^3$ a closed 3-manifold, give a point $x\in Y^3$, for integer $a=1,2,3$, let $\{e_a\}$ be any orthonormal basis of $T_xY$ and let $\{t_a\}$ be elements in the adjoint bundle $ad(\mathfrak{g})$ with the relation that $[t_a,t_b]=\epsilon_{abc}t_c$. For $Y^3\times (0,+\infty)$, we denote by $y$ the coordinate for $(0,+\infty)$.

From [\ref{The Nahm Pole Boundary Condition}] [\ref{Fivebranes and Knots}], we have the following definition of the Nahm pole boundary condition on $Y^3\times (0,+\infty)$.

\begin{definition}
A solution $(A,\phi)$ to the Kapustin-Witten equations $(\ref{KW})$ over $Y^3\times (0,+\infty)$ satisfies the Nahm pole boundary if there exist $\{e_a\}$, $\{t_a\}$ as above such that the Taylor expansion in $y$ coordinate nears $y=0$ will be
$\phi\sim\frac{\sum_{a=1}^3t_ae_a^{\star}}{y}+\phi_0+\ldots,A\sim A_{0}+ya_1+\ldots$.
\end{definition}

\subsection{Nahm Pole Boundary condition over $S^3\times (0,+\infty)$}
In this subsection, we will describe the Nahm pole boundary condition on $S^3$ and show that our solutions in Section 3 satisfy the Nahm pole boundary condition.

Now, we first describe the tangent space of $S^3$. Consider $S^3$ as the unit quaternion, $S^3=\{x=x_1+x_2I+x_3J+x_4K\in \mathbb{H}\vert x_1^2+x_2^2+x_3^2+x_4^2=1\}$, and with the metric induced by Euclidean metric on $\mathbb{R}^4$. Fix a point $x\in S^3$. We can identify the tangent space with vectors on $\mathbb{R}^4$, $T_xS^3=\{v\in \mathbb{R}^4\vert \langle v,x \rangle=0\}$, here we consider $x$ as a vector in $\mathbb{R}^4$.

Define three orthnonormal basis as the following:
\begin{equation}
\begin{split}
e_1&=(-x_2,x_1,-x_4,x_3)\\
e_2&=(-x_3,x_4,x_1,-x_2)\\
e_3&=(-x_4,-x_3,x_2,x_1).
\end{split}
\end{equation}
Obviously, we have $T_xS^3=span\{e_1,e_2,e_3\}$.

Using the induced metric from the Euclidean metric on $\mathbb{R}^4$, we have the dual unit basis
\begin{equation}
\begin{split}
e_1^{\star}&=-x_2dx_1+x_1dx_2-x_4dx_3+x_3dx_4\\
e_2^{\star}&=-x_3dx_1+x_4dx_2+x_1dx_3-x_2dx_4\\
e_3^{\star}&=-x_4dx_1-x_3dx_2+x_2dx_3+x_1dx_4.
\end{split}
\end{equation}

In addition, we have
\begin{equation}
\begin{split}
Im(\bar{x}dx)=&(-x_2dx_1+x_1dx_2-x_4dx_3+x_3dx_4)I\\
&+(-x_3dx_1+x_4dx_2+x_1dx_3-x_2dx_4)J\\
&+(-x_4dx_1-x_3dx_2+x_2dx_3+x_1dx_4)K\\
&=e_1^{\star}I+e_2^{\star}J+e_3^{\star}K.
\end{split}
\end{equation}
Take $t_1=\frac{I}{2}$, $t_1=\frac{J}{2}$, $t_1=\frac{K}{2}$, then we have $[t_a,t_b]=\epsilon_{abc}t_c$.

Therefore, $(0,\frac{Im(\bar{x}dx)}{2})$ can be considered as a leading term of solutions with the Nahm pole boundary condition on $S^3\times (0,+\infty)$.

Now, we will show that the following singular solutions to the Kapustin-Witten equations over $\mathbb{R}^4$ can be viewed as solutions to the Kapustin-Witten equations over $S^3 \times (0,+\infty)$.

As the first equation of ($\ref{KW}$) is conformal invariant, consider the solutions
\begin{equation}
\left\{
\begin{split}
A(x)=&\Im\Big( \frac{3C}{C^2|x|^4+4C|x|^2+1}\;\bar{x}dx\Big)\\
\phi(x)=&\Im\Big( \frac{3C(C|x|^2+1)}{(C^2|x|^4+4C|x|^2+1)(C|x|^2-1)}\;\bar{x}dx\Big).
\end{split}
\label{SolutionNahmpole1}
\right.
\end{equation}
We apply the following conformal transformation
$$\Psi:(0,+\infty)\times S^3\rightarrow \mathbb{R}^4_{\vert x\vert\geq \frac{1}{\sqrt{C}}}$$
$$(y,\omega)\rightarrow (\frac{1}{\sqrt{C}}e^y\omega)$$

The pull back of (\ref{SolutionNahmpole1}) by $\Psi$ gives the following solution on $S^3 \times (0,+\infty)$:
\begin{equation}
\left\{
\begin{split}
A(x)=&\frac{6}{e^{4y}+4e^{2y}+1}\;\sum_{a=1}^3t_ae_a^{\star}\\
\phi(x)=& \frac{6(e^{2y}+1)}{(e^{4y}+4e^{2y}+1)(e^{2y}-1)}\;\sum_{a=1}^3t_ae_a^{\star}.
\end{split}
\right.
\end{equation}

It is easy to see that when $y\rightarrow 0$, $\frac{6(e^{2y}+1)}{(e^{4y}+4e^{2y}+1)(e^{2y}-1)}\sim \frac{1}{y}$ and $y\rightarrow +\infty$, the solution exponentially decays.

From Proposition 3.5, we get another solution to the Kapustin-witten equations ($\ref{KW}$),
\begin{equation}
\left\{
\begin{split}
A(x)=&\Im\Big( \frac{1}{|x|^2}\frac{C^2|x|^4+C|x|^2+1}{C^2|x|^4+4C|x|^2+1}\;\bar{x}dx\Big)\\
\phi(x)=&\Im\Big( \frac{3C(C|x|^2+1)}{(C^2|x|^4+4C|x|^2+1)(C|x|^2-1)}\;\bar{x}dx\Big)
\end{split}
\label{SolutionNahmpole2}
\right.
\end{equation}
Following the same process, we get another solution to the Kapustin-Witten equations with the Nahm pole boundary condition and exponential decay.

Then the pull back of (\ref{SolutionNahmpole2}) using $\Psi$ gives the following solution on $(0,+\infty)\times S^3$:
\begin{equation}
\left\{
\begin{split}
A(x)=&\frac{2}{e^{2y}}\frac{e^{4y}+e^{2y}+1}{e^{4y}+4e^{2y}+1}\;\sum_{a=1}^3t_ae_a^{\star}\\
\phi(x)=& \frac{6(e^{2y}+1)}{(e^{4y}+4e^{2y}+1)(e^{2y}-1)}\;\sum_{a=1}^3t_ae_a^{\star}.
\end{split}
\right.
\end{equation}

Therefore, we have the following theorem:
\begin{theorem}
(a) The following is a solution to the Kapustin-Witten equations over $(0,+\infty)\times S^3$ with the Nahm pole boundary condition and instanton number $+\frac{1}{2}$.
\begin{equation*}
\left\{
\begin{split}
A(x)=& \frac{6}{e^{4y}+4e^{2y}+1}\;\sum_{a=1}^3t_ae_a^{\star}\\
\phi(x)=&\frac{6(e^{2y}+1)}{(e^{4y}+4e^{2y}+1)(e^{2y}-1)}\;\sum_{a=1}^3t_ae_a^{\star}
\end{split}
\right.
\end{equation*}
$$$$
(b) The following is a solution to the Kapustin-Witten equations over $(0,+\infty)\times S^3$ with the Nahm pole boundary and instanton number $-\frac{1}{2}$.
\begin{equation*}
\left\{
\begin{split}
A(x)=& \frac{2}{e^{2y}}\frac{e^{4y}+e^{2y}+1}{e^{4y}+4e^{2y}+1}\;\sum_{a=1}^3t_ae_a^{\star}\\
\phi(x)=&\frac{6(e^{2y}+1)}{(e^{4y}+4e^{2y}+1)(e^{2y}-1)}\;\sum_{a=1}^3t_ae_a^{\star}
\end{split}
\right.
\end{equation*}

\end{theorem}
\proof
The computation of instanton number is along the lines of Prop. 4.3, except we integrate from $1$ to $+\infty$ instead of from $0$ to $+\infty$.
\qed

\section{Acknowledgement}
The author greatly thanks Peter Burton, Anton Kapustin, Edward Witten, Jianfeng Lin, Ciprian Manolescu, Rafe Mazzeo and Yi Ni for their kindness and helpful discussions.

\end{document}